\documentstyle[12pt]{article}
\catcode`\@=11
\@addtoreset{equation}{section}

\catcode`\@=12
\newtheorem{Theorem}{Theorem}[section]

\newtheorem{Proposition}{Proposition}[section]

\title{Proofs On Arnold Conjectures}

\author{Renyi Ma \\
Department of Mathematical Sciences \\
Tsinghua University \\
Beijing, 100084\\
People's Republic of China\\
rma@math.tsinghua.edu.cn}

\date { }

\begin{document}
\textwidth=165mm
\textheight=185mm
\parindent=8mm
\frenchspacing
\maketitle

\begin{abstract}
In this article, we give proofs on the Arnold Lagrangian
intersection conjecture on the cotangent bundles, Arnold-Givental
Lagrangian intersection conjecture, the Arnold fixed point
conjecture.
\end{abstract}

\section{Introduction and results}

 Let $M$ be a smooth manifold and let $T^*M$ denote its
cotangent space. Let $pr_M:T^*M\to M$ denote the natural projection.
The canonical or Liouville $1$-form $\lambda_M$ on $T^*M$ is
$$
\lambda_M(V)=\beta(d(pr_{M})(V))\quad for \   V\in T_\beta(T^*M).
$$
The standard symplectic form on $T^*M$ is the $2$-form
$\omega_M=-d\lambda_M$. If $q=(q_1,\dots,q_n)$ are local coordinates
on $M$ then $(q,p)=(q_1,p_1,\dots,q_n,p_n)$ are local coordinates on
$T^*M$, where $(q,p)$ corresponds to the covector
$$
p_1dq_1+\dots+p_ndq_n\in T^*_q M.
$$
In these local coordinates we have
$$
\lambda_M=-pdq=-\sum_j p_jdq_j \quad \ {and}\quad \omega_M=dq\wedge
dp=\sum_j dq_j\wedge dp_j.
$$
It is easy to see that $(T^*M, d\lambda _M)$ is an exact symplectic
manifold. Let $\phi: T^*M \to T^*M$ be a Hamiltonian
symplectomorphism (see\cite{ag}). The Arnold Lagrange intersection
conjecture(see\cite{arf,ag}) in cotangent bundle is a well-known
conjecture in symplectic geometry. We recall the formulation.
Consider a smooth function $f: M\to R$, we denote by $crit(f)$ and
$crit_m(f)$ the number of critical points of $f$ resp. Morse
function. Let $Crit(M)=\min \{crit(f)\}$ and $Crit_m(M)=\min
\{crit_m(f)\}$ where $f$ runs over all smooth functions $M\to R$
resp. Morse functions. Furthermore, let $Int(M,\phi )$ and
$Int_s(M,\phi )$ denote the number of intersection points of
$\phi(M)$ with $M$ resp. $\phi (M)$ intersects $M$ transversally.
Finally, let
$$
Arn-int(M):=\min_\phi{int(M,\phi)}; \ \ and \ \quad
Arn-int_s(M):=\min_\phi{int_s(M,\phi)}
$$
where $\phi$ runs over all Hamiltonian symplectomorphisms $T^*M \to
T^*M$, resp. $\phi$ runs over all Hamiltonian symplectomorphisms
$T^*M \to T^*M$ such that $\phi (M)$ intersects $M$ transversally.
The Arnold Lagrange intersection conjecture claims that $Arn-int
(M)\geq Crit (M)$ and $Arn-int_s (M)\geq Crit _m(M)$. It is well
known and easy to see that $Arn-int (M) \leq Crit M$ and $Arn-int
_s(M) \leq Crit _m(M)$ . Thus, in fact, the Arnold conjecture claims
the equality $Arn-int(M)=Crit (M)$ and $Arn-int_s(M)=Crit _m(M)$.

\begin{Theorem}
Let $(T^*M,\omega _M)$($\omega _M=d\lambda _M$) be the cotangent
bundles of close manifold $M$. Then,
$$
Arn-int(M):=Crit(M); \ Arn-int_s(M):=Crit_m(M),
$$
i.e., the Arnold Lagrange intersection conjecture in cotangent
bundle holds.
\end{Theorem}
  Theorem1.1 in the stable case was
proved by Hofer\cite{ho}, the other methods provided for example in
\cite{la,eg}, for the complete reference, see \cite{el}.

Now we generalize the above definition to the exact Lagrangian
submanifolds in the exact symplectic manifolds. Let $(V',\omega ')$
be an exact symplectic manifold with exact symplectic form $\omega '
=d\alpha '$. Let $W'\subset V'$ a close submanifold, we call $W'$ an
exact Lagrange submanifold if $\alpha '|W'$ an exact form, i.e.,
$\alpha '|W'=df$. Let $\phi '_t: V'\to V'$ be a Hamiltonian isotopy
with compact support in $V'$ such that $\phi '_0=Id$ and $\phi '
_1=\phi '$ (see\cite{ag}). Let $Int(W',V',\phi ')$ denote the number
of intersection points of $\phi '(W')$ and $W'$ in $V'$. Finally,
let
$$
Arn-int(W',V'):=\min_{\phi '}{int(W',V',\phi ')}
$$
where $\phi '$ runs over all Hamiltonian symplectomorphisms $V'\to
V'$ as above.  Similarly, one defines $ Arn-int_s(W',V'):=\min_{\phi
'}{int_s(W',V',\phi ')}$.

Then the generalized Arnold Lagrange intersection conjecture claims
that $Arn-int (W',V')=Crit (W')$ and $Arn-int_s (W',V')=Crit_m
(W')$.

\begin{Theorem}
Let $(V',\omega ')$ be an exact symplectic manifold with exact
symplectic form $\omega =d\alpha '$. Let $W'\subset V'$ a close
exact Lagrange submanifold. Then,
$$
Arn-int(W',V'):=Crit(W'); \ Arn-int_s(W',V'):=Crit_m(W'),
$$
i.e., the generalized Arnold Lagrange intersection conjecture holds.
\end{Theorem}

      Again we generalize the above definition to the close Lagrangian
submanifold in the general symplectic manifolds. Let $(V',\omega ')$
be a symplectic manifold with symplectic form $\omega '$. Let
$W'\subset V'$ a close submanifold, we call $W'$ a Lagrange
submanifold if $\omega '|W'=0$. Let $\phi '_t: V'\to V'$ be a
Hamiltonian isotopy with compact support in $V'$ such that $\phi
'_0=Id$ and $\phi ' _1=\phi '$ (see\cite{ag}). Let
$Int(W',V',\phi ')$ denote the number of intersection points of
$\phi '(W')$ and $W'$ in $V'$. Finally, let
$$
Arn-int(W',V'):=\min_{\phi '}{int(W',V',\phi ')}
$$
where $\phi '$ runs over all Hamiltonian symplectomorphisms $V'\to
V'$ as above.  Similarly, one defines $ Arn-int_s(W',V'):=\min_{\phi
'}{int_s(W',V',\phi ')}$. Then the Arnold-Givental  Lagrange
intersection conjecture claims that $Arn-int (W',V')=Crit (W')$ and
$Arn-int_s (W',V')=Crit_m (W')$ if $W'$ is symmetric under the
anti-symplectic involution.

\begin{Theorem}
Let $(V',\omega ')$ be a symplectic manifold with symplectic form
$\omega '$. Let $W'\subset V'$ a close Lagrange submanifold such
that $H^1(W')=0$. Then,
$$
Arn-int(W',V'):=Crit(W'); \ Arn-int_s(W',V'):=Crit_m(W'),
$$
\end{Theorem}

\begin{Theorem}
Let $(V',\omega ')$ be a symplectic manifold with symplectic form
$\omega '$. Let $\Phi :(V',\omega ')\to (V',\omega ')$ be a
diffeomorphism such that $\Phi^*\omega '=-\omega '$, $\Phi ^2=Id$.
Let $W'\subset V'$ be the symmetric Lagrange submanifold about
$\Phi$, i.e., $W'=\{x\in V'|\Phi (x)=x\}$. Then,
$$
Arn-int(W',V'):=Crit(W'); \ Arn-int_s(W',V'):=Crit_m(W'),
$$
i.e., the Arnold-Givental Lagrange intersection conjecture holds.
\end{Theorem}

   Now let $(M,\omega)$ be a closed symplectic manifold, and let $\phi: M
\to M$ be a Hamiltonian symplectomorphism (see\cite{ag})).
Furthermore, let $Fix(\phi)$ and $Fix _s(\phi)$ denote the number of
fixed points of general $\phi$ resp transversal $\phi $. Finally,
let
$$
Arn(M,\omega):=\min_\phi{Fix(\phi)} \ and \ \
Arn_s(M,\omega):=\min_\phi{Fix_s(\phi)},
$$
where $\phi$ runs over all Hamiltonian symplectomorphisms $M \to M$
resp. $\phi$ runs over all Hamiltonian symplectomorphisms $M \to M$
such that its fixed points are all nondegenerate. The famous Arnold
conjecture claims that $Arn(M,\omega)=Crit (M)$ and
$Arn_s(M,\omega)=Crit _m(M)$. Let $(V',\omega ')=(M\times M,\omega
\ominus \omega)$, then $L=\{(\sigma ,\sigma)\in M\times M|\sigma \in
M\}$, i.e., its diagonal of the product is Lagrangian submanifold in
$(V',\omega ')$. Moreover, the map $P:(M\times M,\omega \ominus
\omega)\to (M\times M,\omega \ominus \omega),P(x,y)=(y,x),x,y\in M$
is an anti-symplectic involution and $L$ is its fixed points set.
So, by the above theorem, we have

\begin{Theorem}
Let $(M,\omega )$ be a close symplectic manifold with symplectic
form $\omega $. Then,
$$
Arn(M,\omega):=Crit(M) \ and \ \ Arn_s(M,\omega):=Crit _m(M),
$$
i.e., the Arnold fixed point conjecture holds.
\end{Theorem}

\begin{Theorem}
Let $(V',\omega ')$ be a symplectic manifold with symplectic form
$\omega '$. Let $W'\subset V'$ a close Lagrange submanifold such
that $i_*:H_1(W')\to H_1(V')$ is injective. Then,
$$
Arn-int(W',V'):=Crit(W'); \ Arn-int_s(W',V'):=Crit_m(W'),
$$
\end{Theorem}

    The proofs of these theorems does not depend on the hard analysis and 
relies on the analysis of the primitive one form of the symplectic form.

\section{Proof of Theorem 1.1-1.2}

{\bf Proof of Theorem1.1:} Assumption $n\geq 2$.   Consider the
exact isotopy of Lagrange submanifolds in $T^*M$ given by $\varphi
_t,t\in [0,1]$ as a hamilton isotopy of $T^*M$ induced by hamilton vector field 
$X_{h_t}$ with hamilton function $h_t$.
Consider the Lie derivative  $L_{X_{h_t}}\lambda _M$,
by hamilton perturbation, we can assume that
$L_{X_{h_t}}\lambda _M|({\varphi }_t(M)\cap U(M))=
\lambda _M|(\varphi _t(M)\cap U(M))$, 
here $U(M)$ is the small neighbourhood of zero section $M$ in $T^*M$. 
Moreover, $L_{X_{h_t}}\lambda _M=i_{X_{h_t}}d\lambda _M+d 
i_{X_{h_t}}\lambda _M=
dH_t$, here
$H_t$ is defined on $T^*M$. The level sets of
$H_t$ defines a foliation ${\cal {F}}_t$ on
$T^*M$. For $t$ small enough, the Lagrangian submanifold
$\varphi _t(M)$ in $T^*M$ is transversal
to the foliation ${\cal {F}}_t$ except the Lagrangian submanifold
$M$. So, we can perturb the
Lagrangian isotopy ${\varphi }_t$ such that the intersection points sets
$\varphi _t(M)\cap M$ is invariant and
the Lagrangian submanifold
${\varphi }_t(M)$ is transversal
to the foliation ${\cal {F}}_t$ except the intersection points set
$\varphi _t(M)\cap M$.

    This shows that the critical points of $H_1|\varphi _1(M)$ are
in the intersection points set
$\varphi _1(M)\cap M$.

\vskip 3pt

This yields Theorem1.1.

{\bf Proof of Theorem1.2:} Same as the proof of Theorem1.1.

\section{Proof of Theorem1.3-1.6}

\begin{Proposition}
Let $(V',\omega ')$ be a symplectic manifold with symplectic form
$\omega '$. Let $W'\subset V'$ a close Lagrange submanifold. Let
$\varphi '_t: V'\to V'$ be a Hamiltonian isotopy with compact
support in $V'$ such that $\varphi '_0=Id$ and $\varphi ' _1=\varphi
'$. Consider the isotopy of Lagrange submanifolds in $V'$ given by
$F'_t=\varphi _t,t\in [0,1]$ as a $C^{\infty }-$map $F':W'\times
[0,1]\to V'$. We assume that $F':[0,1]\times W'\setminus
\{q_i|i=1,...k\} \to V'$ is an immersion and $F':([0,1]\times
W')\setminus \Sigma \to V'$ is an embedding, here $\Sigma $ is the
Riemann surfaces with boundaries contained in the boundaries  of
$[0,1]\times W'$(see\cite{hi}). Let $L_{n+1}=\cup _{t\in [0,1]}\varphi _t(W')$,
then there exists a neighbourhood $U'(L_{n+1})$ of $L_{n+1}$ in $V'$
and a $1-$form $\alpha '$ on $U'(L_{n+1})$ such that $\omega
'=d\alpha '$.
\end{Proposition}
Proof. With loss of generality, we assume $n\geq 4$. Now we compute
$F'^*\omega $ as in \cite{gro} as
\begin{eqnarray}
F'^*\omega '|([0,1]\times W')&=&\varphi _t^*\omega |W'+ i_{\dot
\varphi }\omega \wedge dt \cr &=&\omega |W'+d_{W'}h_t\wedge dt\cr
&=&d\theta '.
\end{eqnarray}
So, $\omega |F'([0,1]\times W'\setminus (\Sigma \cup
\{q_i|i=1,...k\})$ is exact. Let $\Sigma '\subset U'(L_{n+1})$ is a
$2-$cycle in $W'$, if necessay, by the small perturbation, we can
assume that $\Sigma '\cap (\Sigma \cup \{q_i|i=1,...k\}) =\emptyset
$, so we have $\int _{\Sigma '}\omega '=\int_ {\Sigma '}d\theta
'=0$, so $\omega '|U'(L_{n+1})$ is exact. This proves the
proposition.

\vskip 3pt

{\bf Proof of Theorem1.3.} Since $H^1(W')=0$, by Proposition 3.1 and
Theorem1.2, it is obvious.

\vskip 3pt

{\bf Proof of Theorem1.4.} Since $W'$ is the fixed points set of
$\Phi $ and $\Phi ^*\omega '=-\omega '$. By the proof of
Proposition3.1, there exists a neighbourhood $U'(L_{n+1} \cup \Phi
(L_{n+1}))$ of $L_{n+1}\cup \Phi (L_{n+1})$ in $V'$ and a $1-$form
$\alpha '$ on $U'(L_{n+1} \cup \Phi (L_{n+1}))$ such that $\omega
'=d\alpha '$. Let $\theta '={{1}\over {2}}(\alpha '-\Phi ^*\alpha
')$, then $\omega '=d\theta '$ and $\theta '|W'=0$. By Theorem1.2,
Theorem1.4 follows.

{\bf Proof of Theorem1.6.} By the proof of Proposition3.1, there
exists a neighbourhood $U'(L_{n+1})$ of $L_{n+1}$ in $V'$ and a
$1-$form $\alpha '$ on $U'(L_{n+1})$ such that $\omega '=d\alpha '$.
Then, by the condition that 
 $i_*:H_1(W')\to H_1(V')$ is injective, we can assume that $\alpha
'$ is exact. Then, Theorem1.2 yields
Theorem1.6.

\end{document}